\documentstyle[11pt,psfig]{article}

\input amssym.def
\input amssym.tex

\newcommand{\ga}{\alpha}     
\newcommand{\gb}{\beta}      
\renewcommand{\gg}{\gamma}   
\newcommand{\gd}{\delta}     
  
\newcommand{\gz}{\zeta}      
      
\newcommand{\gth}{\theta}    
      
\newcommand{\gk}{\kappa}  
\newcommand{\gl}{\lambda}    
\newcommand{\gm}{\mu}

\newcommand{\gr}{\rho}       
\newcommand{\gs}{\sigma}     
\newcommand{\gt}{\tau}

\newcommand{\go}{\omega}

\newcommand{\gP}{\Pi}

\newcommand{\ha}{\aleph}

\newcommand{\setof}[2]{{\{\; #1 \; \vert \; #2 \; \} } }
\newcommand{\seq}[1]{{\langle #1 \rangle} }
\newcommand{\card}[1]{{\vert #1 \vert} }
\newcommand{\ot}[1]{\hbox{o.t.($#1$)}}

\renewcommand{\models}{\vDash}

\newcommand{\cf}{{\rm cf}}

\newcommand\pcf{{\rm pcf}}
\newcommand\tcf{{\rm tcf}}

\newcommand{\lra}{\longrightarrow}

\newtheorem{definition}{Definition}
\newtheorem{theorem}{Theorem}
\newenvironment{proof}{\noindent{\bf Proof:}}{\nopagebreak\mbox{}\newline\makebox[\textwidth]{\hfill$\blacklozenge$}\par\bigskip}

\newcommand{\implies}{\Longrightarrow}

\catcode`\@=11

\def\@begintheorem#1#2{
\refstepcounter{subsection}
\rm \trivlist 
\item[\hskip \labelsep{\bf \thesubsection\ #1.}]}

\def\@opargbegintheorem#1#2#3{
\refstepcounter{subsection}
\rm \trivlist
\item[\hskip \labelsep{\bf \thesubsection\ #1\ (#3).}]}

\catcode`\@=12

\newtheorem{lemma}{Lemma}

\newtheorem{corollary}{Corollary}

\newcommand{\FP}{{\Bbb P}}

\newtheorem{remark}{Remark}

\newcommand{\On}{{\it On}}

\newcommand{\Ch}{{\it Ch}}
\newcommand{\Reg}{{\it Reg}}

\newcommand{\Sk}{{\it Sk}}

\newcommand{\CP}{{\cal P}}

\newcommand{\GB}{{\frak B}}

\title{Strong covering without squares}
\author{{\bf Saharon Shelah}\thanks{\ The research was partially supported by
``Basic Research Foundation'' of the Israel Academy of Sciences and
Humanities. Publication 580}
\\
Institute of Mathematics\\
The Hebrew University of Jerusalem\\
91904 Jerusalem, Israel\\
and\\
Department of Mathematics\\
Rutgers University\\
New Brunswick, NJ 08903, USA}

\begin{document}

\maketitle

\baselineskip13.14 truept

\section{Introduction}

  The study of ``covering lemmas'' started with Jensen \cite{Marginalia}
 who proved in 1974--5 that in the absence of $0^\sharp$ there is a certain
 degree of resemblance between $V$ and $L$.
 More precisely, if $0^\sharp$ does not exist then for every set of
 ordinals $X$ there exists a set of ordinals $Y \in L$ such that 
 $X \subseteq Y$ and 
 $V \models \card{Y} = \max \{ \card{X}, \ha_1 \}$. There is no hope of 
 covering countable sets by countable ones in general, because doing
 Namba forcing over $L$ will change the cofinality of $\ha_2^L$
 to $\go$ while preserving $\ha_1$.

 This form of covering has strong implications for the structure
 of $V$. For example Jensen's theorem implies that in the absence
 of $0^\sharp$ the Singular Cardinals Hypothesis holds, and that 
 there is a special $\kappa^+$-Aronszajn tree for every singular
 $\kappa$. So we can conclude that the negations of these 
 statements have substantial consistency strength.

 One subsequent line of development has involved proving covering
 lemmas over larger and larger ``core models'', on the assumption
 of the non-existence of stronger and stronger large cardinals.
 Inevitably these covering lemmas have much more complex
 statements than Jensen's original theorem, the reason being
 that once the core model contains even one measurable
 cardinal we can start to do Prikry forcing.

 This line of research has provided much information
 about the consistency strengths of combinatorial
 hypotheses. 
 For example work of Gitik and Mitchell has determined the exact
 strength of the failure of the Singular Cardinals Hypothesis,
 while work of Schimmerling, Mitchell and Steel has provided
 a very strong lower bound (a Woodin cardinal) for the
 strength of ``$\gk$ is singular and there is no special
 $\gk^+$-Aronszajn tree''. 

 Another line of development involved getting more information
 about the nature of the covering set $Y$. 
 For example given an ordinal $\gl$ and some first-order structure
 ${\cal M} \in V$ with underlying set $\gl$, we can take a set
 $X \subseteq \gl$ and ask whether it can be covered by some
 $Y \in W$ with $\card{Y} \le \max \{ \card{X}, \ha_1 \}$ and
 $Y \prec {\cal M}$. This kind of phenomenon is called
 {\em strong covering} (see Definition \ref{strcov}).

 One approach to proving strong covering theorems is to go
 back to Jensen's proof and to  prove directly that there
 exists an appropriate $Y \in L$. This approach was taken by
 Carlson \cite{Carlson}. Another approach (due to the author)
 is more axiomatic; given $W \subseteq V$ two transitive
 class models of ZFC where $W$ is sufficiently $L$-like
 and for every $X \in V$ there is $Y \supseteq X$
 with $\card{Y} = \max \{ \card{X}, \ha_1 \}$, it is proved
 in \cite[XIII]{ProperForcing}  that a certain form of strong covering
 holds between $V$ and $W$.

 The work in this paper continues that in \cite[XIII]{ProperForcing}
 and \cite{410} 
 (note that a slightly improved version of
\cite[XIII]{ProperForcing} has appeared as
 \cite[VII]{CardinalArithmetic}).
  The idea here is to eliminate as far as possible the structural
 assumptions on $W$. 
 We start by outlining the structure of the paper.
\begin{definition} \label{strcov}
 Let  $W$ be an inner model of ZFC.
 Let $\gk$ be a cardinal in $V$.
\begin{enumerate}
\item {\em $\gk$-covering holds between $V$ and $W$}
 iff for all $X \in V$ with $X \subseteq ON$ and $V \models \card{X} < \gk$,
there exists $Y \in W$ such that $X \subseteq Y \subseteq ON$ and
 $V \models \card{Y} < \gk$.
\item {\em Strong $\gk$-covering holds between $V$ and $W$} iff for every
  structure ${\cal M} \in V$ for some countable first-order language
 whose underlying set is some ordinal $\gl$, and every $X \in V$ with $X \subseteq \gl$
 and $V \models \card{X} < \gk$, there is $Y \in W$ such that $X \subseteq Y \prec M$
 and $V \models \card{Y} < \gk$.
\end{enumerate}
\end{definition}

 In the first section it is proved that if
 $\gk$ is $V$-regular, $\gk^+_V = \gk^+_W$, and we have both $\gk$-covering
 and $\gk^+$-covering between $W$ and $V$, then strong $\gk$-covering holds.
 In fact something rather stronger is proved.
 The assumption that $\gk$-covering holds is reasonable enough, but we can
 hope to weaken the other assumptions.

 In the remainder of the paper we will prove a series of facts about covering
 culminating in two main results; one result says that we can drop the assumption
 of $\gk^+$-covering at the expense of assuming some more absoluteness
 of cardinals and cofinalities between $W$ and $V$, and the other says that
 we can drop the assumption that $\gk^+_W = \gk^+_V$ and weaken the 
 $\gk^+$-covering assumption at the expense of assuming some structural
 facts about $W$ (the existence of certain square sequences).
 Both these results are contained in Theorem \ref{da_big_one}.

 The paper was written up by
 Uri Abraham and James Cummings, and I am grateful for their excellent
 work. I am also grateful to Moti Gitik for asking me about the possibility of
 a theorem like Theorem  \ref{da_big_one}
  after reading \cite{420}.

 The material in this paper represents part of some lectures given
 by the author in Jerusalem in the period May--August 1995.
 The rest of those lectures will appear in \cite{598} and
 so we have retained
 here to some extent the notation and
 terminology used in the lectures.

 In particular the Jerusalem lectures introduced names 
 for some of the important
 hypotheses. For the record, here is a complete list of those names.
 In the body of the paper we will recall these definitions as and
 when we need them.
 $W$ will always be some inner model of ZFC.

\bigskip

\noindent ${\rm (A)}_{\gm,\gl}$: $\ha_0 < \gm = \cf_V(\gm) < \gl \in CARD^W$.

\medskip

\noindent ${\rm (B)}_{\gm,\gl}$: $[\gl]^{<\gm}_W$ is cofinal
 in $[\gl]^{<\gm}$.

\medskip

\noindent ${\rm (B)}_\gm$: ${\rm (B)}_{\gm,\gl}$ holds for every $\gl$.

\medskip

\noindent ${\rm (B)}^-_{\gm, \gl}$: for all $\gd \le \gl$,
 $\cf^W(\gd) < \gm$ iff $\cf(\gd) < \gm$

\medskip

\noindent ${\rm (B)}^*_{\gm, \gl}$: for all $A \in [\gl]^{<\gm}$
 there is $B \in [\gl]^{<\gm^+_V}_W$ such that $A \subseteq B$.

\medskip

\noindent ${\rm (C)}_{\gm,\gl}$: For all $\goth a \in W$,
 if ${\goth a} \subseteq REG^W \cap \gl \setminus \gm$
 then $\max \pcf^W(a) \le \gl$.

\medskip

\noindent ${\rm (D)}_{\gm,\gl}$: In $W$ there is a map
 $\gth \longmapsto \vec C^\gth$ such that
 for each $\gth \in REG^W \cap (\gm, \gl]$,
  $\vec C^\gth = \seq{C^\gth_\ga: \ga < \gth,
 \cf_W(\ga) < \gm}$ is a sequence such that $C^\gth_\ga$ is club in $\ga$,
 $\ot{C^\gth_\ga} < \gm$, and
 $\gb \in \lim(C^\gth_\ga) \implies C^\gth_\ga \cap \gb = C^\gth_\gb$.
 (A {\em square in $\gth$ on points of cofinality less than $\gm$}).

\medskip

\noindent ${\rm (E)}_{\cal P}$: $W \models$ ``${\cal P}\subseteq [\gl]^{<\gm}$,
 $\card{\cal P} = \gl$''. 

\medskip

\noindent ${\rm (F)}_{{\Bbb D}, \gl, \gm}$: For every $x \in [\gl]^\gm$
 there exists
 $\seq{a_i: i < \gm}$ with $a_i \in [\gl]^{<\gm}$ increasing and continuous
 such that $\card{x \cap \bigcup_i a_i} = \gm$,
 $\setof{i < \gm}{a_i \in W} \in {\Bbb D}$.
 Here $\Bbb D$ is a normal filter on
 $\gm$. If $\Bbb D$ is the club filter on $\gm$, we omit it.

\medskip

\noindent ${\rm (G)}_{{\Bbb D}, {\cal P}, \gm}$: for every increasing sequence $\seq{b_i: i < \gm}$ 
 with $b_i \in {\cal P}$ for all $i$, $\setof{i < \gm}{\bigcup_{j < i} b_j \in W} \in {\Bbb D}$.

\medskip

\noindent ${\rm (H)}_{\Bbb D, \gm}$: For all $\ga \in (\gm, \gm^+)$, if $\ga = \bigcup_{i < \gm} a_i$
 with $\vec a$ increasing and continuous 
 then $\setof{i < \gm}{a_i \in W} \in {\Bbb D}$.  
 Here $\Bbb D$ is a normal filter on
 $\gm$. If $\Bbb D$ is the club filter on $\gm$, we omit it.

\medskip

\noindent ${\rm (I)}_{\gm, \gl}$: If $W \models$ ``${\goth a} \subseteq (\gl+1) \setminus \gm^+$,
 $\card{\goth a} < \gm$'' then $\max \pcf^W(\goth a) \le \gl$.

\medskip

\noindent ${\rm (J)}_{\gm,\gl}$: Given a sequence $\seq{{\goth a}_\ga: \ga < \gm}$ such that
 ${\goth a}_\ga \in W$,  ${\goth a}_\ga \subseteq REG^W \cap \gl \setminus \gm^+$,  $\card{{\goth a}_\ga} < \gm$,
 along with $f \in \gP(\bigcup_\ga {\goth a}_\ga)$, there exists $g \in \gP(\bigcup_\ga {\goth a}_\ga)$
 such that $f \le g$ and $g \restriction {\goth a}_\ga \in W$ for all $\ga < \gm$.

\medskip

\noindent ${\rm (K)}_{\gm,\gl}$: $W \models \cf([\gl]^{<\gm}, \subseteq) = \gl$.

\bigskip

   The logical structure of the paper is given by the following picture.

\medskip

\setlength{\unitlength}{0.01in}%

\begin{figure}[h]

\psfig{figure=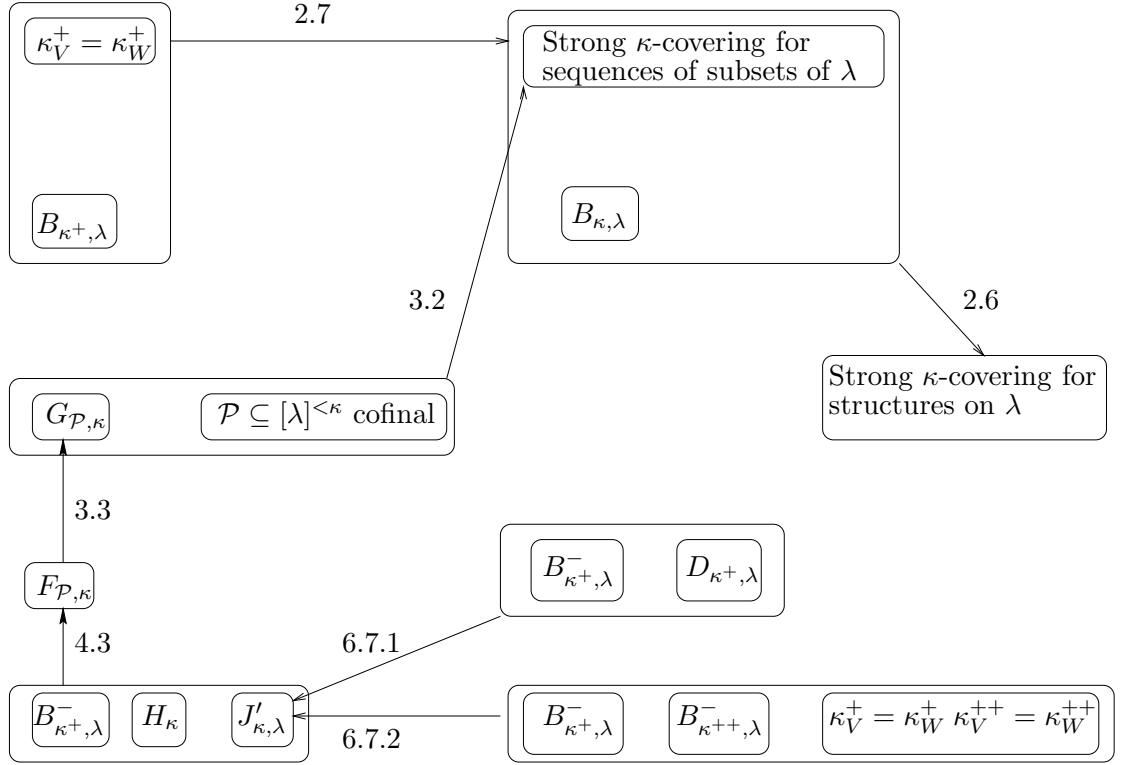,height=4in,rwidth=0in,rheight=0in}

\begin{picture}(600,400)

\put(15,375){$\kappa^+_V=\kappa^+_W$}

\put(150, 390){\ref{27}}

\put(280,375){Strong $\kappa$-covering for}

\put(280,360){sequences of subsets of $\lambda$}

\put(430,200){Strong $\kappa$-covering for}

\put(430,185){structures on $\lambda$}

\put(500,240){\ref{26}}

\put(15,280){$B_{\kappa^+,\lambda}$}

\put(295,285){$B_{\kappa,\lambda}$}

\put(210,240){\ref{Gimpliescovseq}}

\put(20,180){$G_{{\cal P}, \kappa}$}

\put(110,180){${\cal P} \subseteq [\lambda]^{<\kappa}$ cofinal}

\put(35,130){\ref{FimpliesG}}

\put(15,90){$F_{{\cal P}, \kappa}$}

\put(280,100){$B^-_{\kappa^+,\lambda}$}

\put(355,100){$D_{\kappa^+,\lambda}$}

\put(35, 60){\ref{T21}}

\put(175,60){\ref{cor}.\ref{cor1}}

\put(175,10){\ref{cor}.\ref{cor2}}

\put(13,22){$B^-_{\kappa^+,\lambda}$}

\put(70,22){$H_\kappa$}

\put(120,22){$J'_{\kappa,\lambda}$}

\put(280,22){$B^-_{\kappa^+,\lambda}$}

\put(350,22){$B^-_{\kappa^{++},\lambda}$}

\put(430,22){$\kappa^+_V=\kappa^+_W$}

\put(495,22){$\kappa^{++}_V=\kappa^{++}_W$}

\end{picture}

\caption{The structure of the proof}

\end{figure}

\medskip

  We conclude with a few words about notation. If $\gt$ is some set theoretic
 term then ``$\gt_M$'' will mean ``the result of interpreting $\gt$ in the
 model $M$''. When we write ``$M \models \phi(\gt)$'' we mean that 
 ``$\phi$ holds of $\gt_M$ in the sense of the model $M$''. For example
 ``$L\models \cf(\ga) = \go_2$'' means the same thing as
``$\cf_L(\ga) = (\go_2)_L$''.

\section{The painless strong covering theorem}

  In this section we will give a simple proof of a form of strong covering from rather
 strong hypotheses. We begin by discussing the well-known concept of a ``filtration'' which will
 be useful at several points in what follows.

\begin{definition} Let $\gm$ be a regular uncountable cardinal and let $X$ be a set of cardinality
 $\gm$. Then a {\em filtration of $X$\/} is a sequence
 $\vec X = \seq{X_\ga: \ga < \gm}$
 such that
\begin{enumerate}
\item $\card{X_\ga} < \gm$ for all $\ga < \gm$.
\item $\ga < \gb \implies X_\ga \subseteq X_\gb$.
\item For limit $\gl$, $X_\gl = \bigcup_{\ga < \gl} X_\ga$.
\item $X = \bigcup_{\ga < \gm} X_\gm$.
\end{enumerate}
\end{definition}

 The key fact about filtrations is that they are in a sense unique.

\begin{lemma} \label{filtuniq}
 If $\vec X$ and $\vec X'$ are two filtrations of a set $X$ of cardinality $\gm$,
 then $A =_{\rm def}
 \setof{\ga < \gm}{X_\ga = X'_\ga}$ is closed and unbounded in $\gm$.
\end{lemma} 

\begin{proof} For the closure, let $\gl$ be a limit point of $A$,
 and observe that then 
\[
   X_\gl = \bigcup_{\ga < \gl} X_\ga =
  \bigcup_{\ga \in A \cap \gl} X_\ga =
  \bigcup_{\ga \in A \cap \gl} X'_\ga = X'_\gl,
\]
 so $\gl \in A$.
 For unboundedness, fix $\ga_0 < \gm$. Now $X_{\ga_0} \subseteq X = \bigcup_{\ga<\gm} X'_\ga$,
 $\card{X_{\ga_0}} < \gm$ and $\gm$ is regular, so for some $\gb_0< \gm$ we have
 $X_{\ga_0} \subseteq X'_{\gb_0}$. Arguing in a similar vein we may find
 $\ga_0 < \gb_0 < \ga_1 <\gb_1 \ldots$ such that
 $X_{\ga_0} \subseteq X'_{\gb_0} \subseteq X_{\ga_1} \subseteq X'_{\gb_1} \ldots$.
 If now $\gl = \sup_n \ga_n = \sup_n \gb_n$, then 
 $X_\gl = \bigcup_n X_{\ga_n} = \bigcup_n X'_{\gb_n} = X'_\gl$ so that
 $\gl \in A$ and we have proved that $A$ is unbounded.
\end{proof}

 Recall that the hypothesis
 ${\rm (B)}_\gm$ says that if $X \in V$, $X \subseteq ON$
 and $V \models \card{X} < \gm$
 then there is $Y \in W$, $X \subseteq Y \subseteq ON$
 with $V \models \card{Y} < \gm$. This is the hypothesis which we called
 ``$\gm$-covering'' in the introduction.

\begin{theorem} \label{painless}
 Let $\gk$ be regular and uncountable in $V$. Suppose
 that ${\rm (B)}_\gk$
 and ${\rm (B)}_{\gk^+}$ both hold and that $\gk^+_V = \gk^+_W$.
 Let ${\cal M} \in V$ be
 a structure for some countable first order language, whose
 underlying set is an
 ordinal $\gl \ge \gk$.
 Let $z$ be a set in $[\gl]^{<\gk}$.
 Then there exists an increasing and
 continuous sequence of sets $\seq{z_\ga: \ga < \gk}$ 
 such that $z \subseteq z_0$,
$\card{z_\ga} < \gk$, $z_\ga \prec {\cal M}$, and $z_\ga \in W$.
 In particular, strong $\gk$-covering holds between $V$ and $W$.
\end{theorem}

\begin{proof} We build an increasing and
 continuous sequence $\seq{x_\ga: \ga < \gk}$
 of subsets of $\gl$, where each $x_\ga$ has size less than $\gk$.

\begin{itemize}

\item In case $\ga = 0$,  $x_0 = z$.

\item In case  $\ga = 2 \gb +1$,  $x_\ga$ is some set such that
 $x_{2\gb} \subsetneq x_\ga \subseteq \gl$,
 $x_\ga \in W$, $\card{x_\ga} < \gk$. Such a set
 exists because we are assuming ${\rm (B)}_\gk$.

\item In case  $\ga = 2 \gb +2$, 
 $x_\ga$ is some set such that $x_{2\gb+1} \subsetneq x_\ga \subseteq \gl$, 
 $x_\ga \prec {\cal M}$, and $\card{x_\ga} < \gk$.
 Such a set exists because $\cal M$ has a countable set of Skolem functions
 and $\gk$ is uncountable.

\item In case  $\ga$ is limit, $x_\ga = \bigcup_{\gb < \ga} x_\gb$.
 $\card{x_\ga} < \gk$ because $\gk$ is regular.
 Observe that $x_\ga \prec {\cal M}$,
 because $x_\ga$ is the union of the increasing chain of
 substructures $\seq{x_{2\gb+2}: 2\gb+2<\ga}$.

\end{itemize}

  Now let $x = \bigcup_{\ga < \gk} x_\ga$, so that $\card{x} = \gk$.
 Applying the hypothesis ${\rm (B)}_{\gk^+}$, we may
 find $y \subseteq \gl$ such that $y \in W$, $x \subseteq y$ and $V \models \card{y} < \gk^+$.
 Since we know that $\gk^+_V = \gk^+_W$, we see that in fact $W \models \card{y} = \gk$.
 Working in $W$, let us fix a filtration $\seq{y_\ga: \ga < \gk}$ of $y$.

 Now we observe that we have two filtrations of the set $x$, because
 $x = \bigcup_\ga x_\ga = \bigcup_\ga (x \cap y_\ga)$. By Lemma \ref{filtuniq}
 we can find a club $C$ of limit ordinals in $\gk$ such that $\ga \in C \implies x_\ga = x \cap y_\ga$.
 Now $x_\ga \subseteq x_{\ga+1} \subseteq x$ so that for $\ga \in C$
\[
    x_\ga = x_\ga \cap x_{\ga+1} = x \cap y_\ga \cap x_{\ga+1} = y_\ga \cap x_{\ga+1},
\]
 and hence $x_\ga \in W$ since both $x_{\ga+1}$ and $y_\ga$ are in $W$.
 If we enumerate $C$ in increasing order as $\seq{\gg_\ga: \ga < \gk}$ and set $z_\ga = x_{\gg_\ga}$
 then $\seq{z_\ga: \ga < \gk}$ is continuous and increasing, $z_\ga \prec M$ and $z_\ga \in W$
 for all $\ga < \gk$. 
\end{proof}

\begin{remark}
 The assumption ${\rm (B)}_\gk$ in the theorem is redundant, because
 actually the assumptions ${\rm (B)}_{\gk^+}$ and $\gk^+_V = \gk^+_W$
 imply  ${\rm (B)}_{\gk}$. For if $\card{X} < \gk$ we may cover it by $Y \in W$
 such that $\card{Y}_V < \gk^+$, but then by the agreement between cardinals
 we must have $\card{Y}_W \le \gk$. Now in $W$ we write $Y = \bigcup_{i<\gk} Y_i$
 with $\card{Y_i} < \gk$, and since $\gk$ is regular in $V$ there exists $i$
 such that $X \subseteq Y_i$.
\end{remark}

 We can analyse the proof of Theorem \ref{painless} into two steps, an analysis which
 gives some motivation for the work of later sections.

\begin{definition}
Let $\gk$ be regular, let $W$ be an inner model.
\begin{enumerate}
\item A set $X \in V$ with $\card{X} = \gk$ and $X \subseteq W$ is {\em $W$-filtered\/} iff
 there is a filtration $\vec X$ of $X$ such that $X_\ga \in W$
 for a closed unbounded set of $\ga < \gk$. 
\item {\em Strong $\gk$-covering for sequences holds between $V$ and $W$ for subsets of $\gl$\/}
 iff whenever
 $\seq{a_i: i < \gk}$ is an increasing (but not necessarily continuous)
 sequence of subsets of $\gl$ with $a_i \in W$ and $\card{a_i} < \gk$,
 then $\bigcup_i a_i$ is $W$-filtered.
\item {\em Strong $\gk$-covering for sequences holds between $V$ and $W$\/} iff 
 for all $\gl$ strong $\gk$-covering for sequences holds between $V$ and $W$ for subsets of $\gl$.
\end{enumerate}
\end{definition}

 The proof of Theorem \ref{painless} can be broken into two lemmas.

\begin{lemma} \label{26}
 Let ${\rm (B)}_\gk$ hold, let $\cal M$ be a structure
 for a countable first-order language with underlying set $\gl$,
 and let $z \subseteq \gl$ with $\card{z} < \gk$. Then 
 $z \subseteq \bigcup_{i<\gk} a_i$ where $\card{a_i} < \gk$,  $a_i \in W$ for
 $i$ odd, $a_i \prec M$ for $i$ even. 
\end{lemma}

\begin{lemma} \label{27}
 If ${\rm (B)}_{\gk^+}$ and $\gk^+_V = \gk^+_W$ then
 strong $\gk$-covering for sequences holds.
\end{lemma}
 It follows from the second lemma that for many limit $i$ we have $a_i \in W$
 and hence we have covered $z$ by a substructure lying in $W$.

\begin{remark} Notice that by Lemma \ref{filtuniq}, if $X$ is $W$-filtered
 then for any filtration $\vec X$ we have $X_\ga \in W$ for a club of $\ga$.
 Another equivalent definition of ``$W$-filtered'' would demand that $X$ has a filtration 
 such that $X_\ga \in W$ for every $\ga$.
\end{remark}

\begin{remark} Strong $\gk$-covering for sequences does not imply $\gk$-covering.
 For example if $L[G]$ is the generic extension of $L$ for Namba forcing then
 strong $\go_1$-covering for sequences holds (by an easy argument using the fact
 that Namba forcing adds no reals to a model of CH) but $\go_1$-covering certainly
 fails.
\end{remark}

\section{${\rm (F)}$ implies ${\rm (G)}$}

  We now begin to show how to weaken the hypotheses of Theorem \ref{painless}
 from the last section.
Throughout this section $\gk$ denotes a regular uncountable cardinal (in
$V$), and $\gl$ a $W$-cardinal with $\gl > \gk$.

 The next lemma represents a variation on the main idea in
 Theorem \ref{painless}. Here we are covering by a set of size
 $\gk$ which is $W$-filtered, rather than actually lying in $W$.

\begin{lemma} \label{inclfilt}
 Suppose that $b= \bigcup_{i< \gk} b_i$ where $b_i \in W$,
$\card{b_i} < \gk$, $b_i \subseteq b_j$ for $i < j < \gk$. If $b$ is
{\em included} in a $W$-filtered set of size $\gk$, then $b$ itself is $W$-filtered.
\end{lemma}

\begin{proof} Suppose that $b \subseteq a$ where $a$ is $W$-filtered, as
evidenced by an increasing and continuous sequence
 $\seq{ a_i : i < \gk}$ such that $a=\bigcup_{i < \gk} a_i$,
 $a_i \in W$, $\card{a_i} < \gk$. Then the sets $a_i \cap b$ form one
 filtration of $b$, while
$\seq{\bigcup_{ i < \gd} b_i : \gd < \gk}$ is another.
By Lemma \ref{filtuniq},
\[
 D= \setof{\gd < \gk}{a_\gd \cap b= \bigcup_{\ga <\gd} b_\ga}
\]
contains a club set. We will prove for every $\gd \in D$ that $a_\gd
\cap b \in W$.
 Since the sequence of $b_\ga$'s is
increasing, $\bigcup_{\ga < \gd} b_\ga \subseteq b_\gd$. Hence since
$ a_\gd \cap b= \bigcup_{\ga < \gd} b_\ga$,
  $a_\gd \cap b = a_\gd \cap b_\gd$ 
as well.  Hence, if $\gd \in  D$, then $a_\gd \cap b  \in W$ 
(since $a_\gd, b_\gd \in W$). 
\end{proof}

Recall that  ${\rm (F)}_{\gl,\gk}$ denotes the following statement:

\medskip

\noindent ${\rm (F)}_{\gl,\gk}$: for all $x \in [\gl]^\gk$ there exists $A$ such that $\card{A} = \gk$,
 $A$ is $W$-filtered and $\card{x \cap A} = \gk$.

\medskip

 Let $\CP$ be such that  $W \models \hbox{``$\CP \subseteq [\gl ] ^{< \gk}$, $\card{\CP} = \gl$''}$.
 In applications  of ${\rm (F)}_{\gl,\gk}$ 
we will replace $\gl$ with
$\CP$, and then  ${\rm (F)}_{\gl,\gk}$ becomes the following sentence

\medskip

\noindent ${\rm (F)}_{\CP,\gk}$:
  for all $x \in [\CP]^\gk$ there exists $A$ such that $\card{A} = \gk$,
 $A$ is $W$-filtered and $\card{x \cap A} = \gk$.

\medskip

Recall also that we defined

\medskip

\noindent  ${\rm (G)}_{\CP, \gk}$:
for every $\subseteq$-increasing sequence of sets  
$\seq{b_\ga : \ga < \gk}$ with $b_\ga \in \CP$ for all $\ga$,
 the set $\setof{\gd < \gk}{\bigcup_{\ga < \gd} b_\ga \in W}$
 contains a club set.

\medskip

Equivalently, ${\rm (G)}_{\CP, \gk}$ says that for every $\subseteq$-increasing
sequence $\langle b_\ga : \ga < \gk \rangle$ with $b_\ga \in \CP$,
if $b=\bigcup _{\ga < \gk} b_\ga$ then $b$ is $W$-filtered.
 This is a weakening of the statement of strong $\gk$-covering for sequences,
 where we restrict to the class of increasing sequences from $\CP$.

\begin{lemma} \label{Gimpliescovseq} 
If $(\CP,\subseteq)$ is cofinal in
$(\gl^{< \gk}, \subseteq)$,
${\rm (G)}_{\CP, \gk}$ implies the strong $\gk$-covering for sequences of subsets of
$\gl$.
\end{lemma}

\begin{proof}
  Let $\seq{a_i : i < \gk}$ be an
increasing sequence of subsets of $\gl$ such that $a_i \in W$ and
 $\card{a_i} < \gk$. It is required to prove that $a=\bigcup_{i< \gk} a_i$ is
$W$-filtered. By Lemma \ref{inclfilt} it suffices to prove that $a$ is
included in a $W$-filtered set, and for this we define inductively an
increasing sequence $b_i \in \CP$ for $i < \gk$ such that $a_i \subseteq
b_i$. Now by ${\rm (G)}_{\CP, \gk}$ the set $b=\bigcup_{i< \gk} b_i$ is
$W$-filtered,  and $a \subseteq b$,  hence $a$ is $W$-filtered.
\end{proof}

\begin{theorem} \label{FimpliesG}
${\rm (F)}_{\CP, \gk}$ implies ${\rm (G)}_{\CP, \gk}$.
\end{theorem}

\begin{proof}
 Let $\seq{b_\ga : \ga < \gk}$ be an increasing
sequence of members of $\CP$, and set $b=\bigcup_{\ga < \gk} b_\ga$.
  Apply ${\rm (F)}_{\CP, \gk}$ to the set
 $x=\setof{b_\ga}{\ga < \gk} \in [ \CP]^\gk$,
 and obtain a $W$-filtered set $A \subseteq \CP$ such that
 $\card{x \cap A} = \gk$.
 Write $A= \setof{a_i}{i<\gk}$, and define $a=\bigcup A$. We will show
that $b \subseteq a$ and $a$ is $W$-filtered. By Lemma \ref{inclfilt} this implies that $b$ is
$W$-filtered, as required.

Since $x \cap A$ contains unboundedly many $b_\ga$'s, 
$b \subseteq a = \bigcup A$.
  Since $A$ is
$W$-filtered,
 $\setof{\gd < \gk}{\setof{a_i}{i < \gd }\in W}$
contains a club set, and hence the larger set
\[
   C= \setof{\gd < \gk}{\bigcup_{i<\gd} a_i \in W}
\]
contains a club set. This gives a witness that $a$ is $W$-filtered.
\end{proof}

\begin{remark} \label{Fiscovering}
 If $(\CP,\subseteq)$ is cofinal in
$(\gl^{< \gk}, \subseteq)$, Theorem \ref{FimpliesG} and
 Lemma \ref{Gimpliescovseq} 
 imply that the strong
covering property for sequences of subsets of $\gl$ follows from
${\rm (F)}_{\CP, \gk}$.  The theorem also shows in this case that ${\rm (F)}_{\CP, \gk}$
implies that every $x \in [\gl]^\gk$ is {\em contained} in a
$W$-filtered set.
\end{remark}

\begin{remark}
 We could have replaced the club filter on $\gk$ by any
filter ${\Bbb D}$ over $\gk$ containing the club sets.  Then a
$W$-$\Bbb D$-filtered set would be a set $A$ with a decomposition 
$A=\bigcup_{i<\gk} a_i$ such that
$\{ \gd < \gk : \bigcup_{i<\gd} a_i \in W\} \in {\Bbb D}$.  The
statement ${\rm (F)}_\gl$ would be replaced by (the weaker) ${\rm (F)}_{\Bbb D}$, and
correspondingly ${\rm (G)}_\gl$ would be replaced by ${\rm (G)}_{\Bbb D}$.  The
proof that ${\rm (F)}_{\Bbb D} \implies {\rm (G)}_{\Bbb D}$ would be almost the same if ``$X$
is a club'' is replaced by $X \in {\Bbb D}$.
 Note: working with an arbitrary normal filter $\Bbb D$ on $\gk$ does not materially
 change the proof, but gives a much
 weaker assumption.
\end{remark}

\section{ How to derive ${\rm (F)}$}
We saw in Remark \ref{Fiscovering} that ${\rm (F)}$ is a form of strong covering.
In this section we see how to derive it from some other putatively
 weaker assumptions.

We start by listing the necessary assumptions.

\medskip

\noindent ${\rm (H)}_\gk$: every ordinal $\ga$ in the interval $(\gk, \gk^+)$
 is $W$-filtered.

\medskip

\begin{remark}  Notice that $H_\gk$ is a weakening of the assumption that
 $\gk^+_V = \gk^+_W$. It is properly weaker, as can be seen by
 starting with a model $W$ and $\gk < \gl$ with $\gk$, $\gl$ both
 regular in $W$, and doing the Levy collapse $Coll(\gk, <\gl)$. 
\end{remark}

\medskip

 Recall that a function $g:
A \rightarrow \On$ from a set of ordinals $A$  of size $\gk$ into the
 ordinals is
$W$-filtered iff there exists $\vec a$ a filtration of $A$
  such that for all $i$ both $a_i$
 and $g \restriction a_i$ are in $W$.

\medskip

\noindent ${\rm (J)'}_{\gk, \gl}$ :
For any ${\frak a} \subseteq \Reg^W \cap ( \gl \setminus
\gk^+)$ of cardinality $\gk$, if ${\frak a}$ is $W$-filtered, then for
every $f\in \Pi {\frak a}$ there exists $g \in \Pi {\frak a}$ such that $f \leq g$ and
$g$ is $W$-filtered.

\medskip

Here $f \in \Pi {\frak a}$ means that $f(\ga) \in \ga$ for every $\ga
\in {\frak a}$, and $f \leq g$ means that $\forall \ga \in {\frak a} \;  (f(\ga) \leq
g(\ga))$. $J'$ is a weaker version of the property $J$ defined
 in the introduction.

\medskip

\noindent ${\rm (B)}^-_{\gk^+,\gl}$:  For every ordinal $\ga \in \gl$, if
 $\cf (\ga)< \gk$, then $\cf_W(\ga) < \gk^+_V$.

\medskip
 
\begin{remark} The hypotheses ${\rm (B)}^-_{\gk^+,\gl}$
 and  ${\rm (J)'}_{\gk, \gl}$ are both consequences of
 ${\rm (B)}_{\gk^+, \gl}$.
\end{remark}

\begin{theorem} \label{T21}
 If	${\rm (B)}^-_{\gk^+,\gl}$,  ${\rm (J)'}_{\gk, \gl}$ and ${\rm (H)}_\gk$
 then
 ${\rm (F)}_{\gl, \gk}$.
\end{theorem}

\begin{proof}  Let $x \in [\gl]^\gk$ be
given;  we shall find a $W$-filtered set $A$ such that $\mid x \cap
A \mid = \gk$.  For this we define by induction
$W$-filtered sets $N_n$ for each $n < \go$, such that $x \subseteq \bigcup_{n < \omega}
N_n$; and then necessarily some $N_i$ satisfies $\mid x \cap N_i \mid =
\gk$ and is thus as required.

First we fix a rich enough first-order structure $\GB\in W$, with universe $\gl$
 and finitely many functions and relations.
 The meaning of ``rich enough'' will become clear
 during the proof, where we will list a finite set of functions which
 should appear among the functions of $\GB$. 

Then we construct elementary submodels $N_n,M_n \prec \GB$ of size $\gk$
such that:
\begin{enumerate}
\item $N_n$ is $W$-filtered, and $\gk \subseteq N_n\subseteq N_{n+1}$.
\item $x \subseteq M_0,\ M_n \subseteq M_{n+1}$ ($x$ is the given set).
\item $N_n \subseteq M_n$.
\item For ${\frak a}_n= \Reg^W \cap (\gl \setminus \gk^+) \cap N_n$, if we define
$f_n=\Ch_{{\frak a}_n}^{M_n}$ and $g_n=\Ch_{{\frak a}_n}^{N_{n+1}}$, then  $f_n \leq g_n$.
\end{enumerate}

Recall that $\Ch_A^M(\theta)=\sup (M\cap\theta)$ for $\theta \in A$.  Thus
$f_n \leq g_n$ means simply that $\sup (M_n \cap \theta) \leq \sup
(N_{n+1}\cap \theta)$, for all $W$-regular $\theta \geq \gk^+$ in $N_n$.

The construction begins by setting $N_0=\Sk(\gk)$, where $\Sk(X)$ is the
Skolem closure of $X \subseteq \gl$ in $\GB$. This maintains the
 induction hypotheses because of the following fact.

\begin{lemma}
If $X$ is $W$-filtered, then  so is $\Sk(X)$.
\end{lemma}

\begin{proof}
Easy, using the fact that $\GB\in W$.
\end{proof}

In the second step, $M_0$ is defined to be any elementary submodel of
$\GB$ of size $\gk$ which contains $x$ and $N_0$.  For example,
$M_0=\Sk(N_0 \cup x)$ will do.

Suppose now that $N_n$ and $M_n$ have been defined.  We first define
$N_{n+1}$ as follows.  Set ${\frak a}_n=\Reg^W \cap(\gl \setminus \gk^+)
\cap N_n$.  For every $\theta \in {\frak a}_n$, as ${\rm (B)}^-_{\gk^+,\gl}$ holds
and as $\theta \geq \gk^+$ is regular in $W$, $\cf_V(\theta) \geq
\gk^+$.  Hence $f_n(\theta)=\sup(M_n \cap \theta) < \theta$, and
the function $f_n$ thus defined is in $\Pi {\frak a}_n$.  So ${\rm (J)'}_\gl$ can
be applied to ${\frak a}_n,f_n$ and there exists $g_n \in \Pi {\frak a}_n$ which is
$W$-filtered and such that $f_n \leq g_n$. It follows that the set $\{
g_n(\theta) : \theta \in {\frak a}_n \}$ is also $W$-filtered.  Let
$\gamma_n=\sup(M_n \cap \gk^+)$, then $\gamma_n$ is $W$-filtered, because 
 we are assuming ${\rm (H)}_\gk$.  Let $\sigma_n=\sup M_n$ if $\sup M_n < \gl$, and
$\sigma_n=0$ otherwise. Define
\[
  N_{n+1}=\Sk(N_n \cup \{ g_n(\theta) : \theta \in {\frak a}_n \} \cup \gamma_n
                       \cup \{ \sigma_n \})
\]
$N_{n+1}$ is $W$-filtered.
Finally, define
\[
    M_{n+1}=\Sk(N_{n+1} \cup M_n).
\]
 This ends the inductive definition of the models, and we  prove
now that
$\bigcup_{n<\omega} N_n= \bigcup_{n< \omega}M_n$.
This obviously implies $x \subseteq \bigcup_{n<\omega}M_n$ and thus ends
the proof.

Let $N=\bigcup_{n<\omega} N_n$, $M= \bigcup_{n< \omega}M_n$.  The
following properties easily follow from our construction.
\begin{enumerate}
\item $N \prec M$. (Since $N_n \prec M_n$.)
\item $N \cap \gk^+=M\cap \gk^+=\sup \{ \gamma_n : n < \omega \}$.
\item $\sup N = \sup M \leq \gl$.
\item For every $\theta \in N \cap \Reg^W \setminus \gk^+$,
$\sup(N\cap \theta)= \sup(M \cap \theta)$.  Indeed $\sup(N\cap
\theta)\leq \sup(M\cap \theta)$ follows from $N \subseteq M$, and the
other direction of the inequality follows from $\sup (M_n \cap \theta )<
\sup (N_{n+1} \cap \theta)$.
\end{enumerate}

 We claim that (1), (2), (3), (4) imply that $M=N$. To see this,
assume otherwise and then $M \setminus N \neq \emptyset$.
  Let $\ga=\min(M \setminus N)$.  By (3) $\beta=\min(N \setminus \ga)$
 is defined.  We shall derive a contradiction in the following
complete analysis of cases for $\beta$. 
\begin{description}

\item[Case:] $\beta < \gk^+$.  Then $\ga < \gk^+$ and this
contradicts (2).

\item[Case:] $\beta$ {\bf is a successor ordinal}.  Let $\tau$ be such
that $\tau+1 =\beta$.  Then $\tau \in N$, if $\GB$ is closed under
predecessors. We will therefore demand that the predecessor function
 appear among the functions of $\GB$.

Now $\ga <\beta$ implies $\ga \leq \tau$, and yet $\ga = \tau$
is not possible since $\ga \not \in N$; hence $\ga < \tau < \beta$
contradicts the minimality of $\beta$.

\item[Case:] $\beta$ {\bf is singular in} $W$.  Let $\tau=\cf^W(\beta)$.
Then $\tau \in N$ as long as the $W$-cofinality function is in $\GB$.
 We therefore demand that $\cf_W$ is among the functions of $W$.

  Now as
$\tau < \beta$ ($\beta$ being singular) $\tau <\ga$ follows. 
We will demand that $\GB$ contains a function ${\rm COF}$ which
assigns to every $W$-singular cardinal $\gz$ a cofinal sequence ${\rm COF}(\gz) \in W$ of 
length $\cf_W(\gz)$. ${\rm COF}(\beta): \tau \lra \beta$ is cofinal in $\beta$ and
so, as $\ga \in M$, there is $\zeta < \tau$ in $M$ such that
$\ga'={\rm COF}(\beta)(\zeta) \geq \ga$.  Since $N \cap \ga=M \cap\ga$
by the minimality of $\ga$, $\zeta \in N$. Hence $\ga' \in N$ is
an ordinal in $[\ga, \beta)$, in contradiction to the minimality of
$\beta$.

\item[Case:]
 $\beta$ {\bf is regular in} $W$.  This is the last case.  We
know by (4) that $\sup(N\cap \beta)=\sup(M \cap \beta)$.  Hence
$\sup(N\cap \beta) \geq \ga$ which is again a contradiction to
$\beta=\min(N\setminus \ga)$.

\end{description}

This concludes the proof of the Theorem.
\end{proof}

  We end this section by considering a variation 
 ${\rm (J)}^{''}_\gl$ which seems weaker than ${\rm (J)'}_\gl$
but suffices for Theorem \ref{T21}

\medskip

\noindent ${\rm (J)}^{''}_\gl$:
If $A \subseteq \Reg^W \cap (\gl \setminus \gk^+)$ of cardinality
$\gk$ is $W$-filtered and $f \in \Pi A$, then there exists $i^* <
\gk$ and a collection $\{ f_i : i < i ^* \}$, $f_i \in \Pi A$, such that
$f \leq \sup \{ f_i : i < i^* \}$ and each $f_i$ is $W$-filtered.

\medskip

\begin{theorem}
Theorem \ref{T21} still holds if ${\rm (J)}^{''}_\gl$ replaces ${\rm (J)'}_\gl$.
\end{theorem}

\begin{proof}

Given $i^* < \gk$, let us say that $A$ is {\em $(i^*,W)$-filtered\/} iff $A$ is a union of $i^*$ sets
each of which is $W$-filtered.

We will define $N_n$, $M_n$ as before, but require now that $N_n$ is
$(i^*_n,W)$-filtered for some $i^*_n < \gk$. The construction is very
similar, but in defining $N_{n+1}$ it is not a single function that is
added to $N_n$, but rather $i^*_n$ functions, each $W$-filtered. Finally
$N=M$ as before, and for some $n$ we have $\card{N_n \cap x} = \gk$.

As $N_n$ is $(i^*_n,W)$-filtered, we may find a $W$-filtered set $X
\subseteq N_n$ such that $\card{X \cap x}= \gk$ as required.
\end{proof}

\section{The pcf induction lemma}

   The results in the last section indicate the usefulness of the
 hypothesis ${\rm (J)'}$. In this section we prove a crucial lemma,
 which we will exploit in the final section of the paper to prove that
 two rather different sets of assumptions will each imply ${\rm (J)'}$.

 \begin{lemma} \label{soft}
 Let $W \subseteq V$ be two transitive class models of
 set theory. In $V$ let $\gk$ and $\gs$ be regular, with $\gk < \gs$.
 Assume that
\begin{enumerate}
\item $N \prec (H_\gs, \in, <^*, W \cap H_\gs)$ where $<^*$ is a
 wellordering of $H_\gs$.
\item $\card{N} = \gk$, $\gk \subseteq N$.
\item If $\gt \in N \cap REG^W$ with $\gt \ge \gk^+_V$, then
\begin{enumerate}
\item  $\sup(N \cap \gt) < \gt$. \label{3a}
\item There exists $C$ an unbounded subset of $N \cap \gt$
 with $C \in W$. \label{3b}
\end{enumerate}
\item \label{4}
 For all ${\frak a} \subseteq REG^W \setminus \gk^+_V$ with ${\frak a} \in N \cap W$,
 $\card{{\frak a}} < \gk$
 and $\max(\pcf_W({\frak a})) < \gs$, if $h \in \Pi {\frak a}$ and $h(\gr) < \sup(N \cap \gr)$
 for all $\gr \in {\frak a}$, then there is $g \in \Pi {\frak a} \cap N \cap W$ such that
 $h(\gr) < g(\gr)$ for all $\gr \in {\frak a}$.
\end{enumerate}

   Then for all ${\frak a} \subseteq REG^W \setminus \gk^+_V$ with ${\frak a} \in N \cap W$,
  $\max(\pcf_W(a)) < \gs$ and $\card{{\frak a}} < \gk$, if we define 
 $ch^N_{\frak a}: \gr \in {\frak a} \longmapsto \sup(N \cap \gr)$
 then $ch^N_{\frak a} \in W$.
\end{lemma}

\begin{proof} Fix a structure $N$ with these properties. 
  We prove the lemma by induction on $\gth = \max(\pcf_W({\frak a}))$. We start the
 induction by observing that if $\gth = \gk^+_V$ then necessarily
 ${\frak a} = \{ \gk^+_V \}$
 (because ${\frak a} \subseteq \pcf({\frak a})$ and $\min({\frak a}) = \min(\pcf({\frak a}))$)
 and the result is trivial.

 For the general case, start by observing that since ${\frak a} \in N$ and
 $\gth < \gs$ we have $\gth \in N$. Notice also that by \ref{3a} we have
 $ch^N_{\frak a} \in \Pi {\frak a}$.

 Fix a sequence $\seq{f_\gg: \gg < \gth} \in N \cap W$ such that
\[
   W \models \hbox{``$\vec f$ is cofinal in $\Pi {\frak a}/J_{<\gth}$''}.
\]
 Such a sequence $\vec f$ exists in $W$ by
 a basic fact from pcf theory
 (see Theorem \ref{bop} from the Appendix)
 and since ${\frak a}, \gth \in N$ we may assume that $\vec f$ lies in $N$.

 $\gth \in N$ and $\gth$ is regular in $W$, so  by \ref{3b} we may fix some
 $C$ unbounded in $N \cap \gth$ with $C \in W$.
 Let us define $f(\gr) = \sup_{\gg \in C} f_\gg(\gr)$ for $\gr \in {\frak a}$. Then
 \begin{enumerate}
\item $f \in W$ because $f$ is defined from $\vec f$, ${\frak a}$  and $C$, which
 are all in $W$.
\item For each $\gg \in C$ and $\gr \in {\frak a}$, $f_\gg(\gr) \in N \cap \gr$
 because $C, {\frak a} \subseteq N$ and $\vec f \in N$. Therefore
 $f(\gr) \le ch^N_{\frak a}(\gr) < \gr$.
\end{enumerate}

 Let ${\frak b} =_{\rm def} \setof{\gr \in {\frak a}}{f(\gr) < ch^N_{\frak a}(\gr)}$,
 and define $h$ by setting $h \restriction {\frak b} = f \restriction {\frak b}$ and
 $h \restriction {\frak b}^c$ to be the zero function.  
  By \ref{4} we may find $g \in\Pi {\frak a} \cap N \cap W$
 such that $h < g$, so that $\forall \gr \in {\frak b} \; f(\gr) < g(\gr)$. Using
 elementarity there is $\gg \in N \cap \gth$ such that
 $g <_{J_{<\gth}} f_\gg$, and since $C$ is unbounded in $N \cap \gth$
 we may assume that $\gg \in C$. By the construction of $f$,
 $f_\gg \le f$.  

 Now let ${\frak b}^* = \setof{\gr \in {\frak a}}{f_{\gg}(\gr) \le g(\gr)}$. Then
 ${\frak b}^* \in N \cap W$ because the parameters ${\frak a}$, $f_\gg$ and $g$ lie
 in $N \cap W$.  ${\frak b}^* \in J_{<\gth}$, because $g <_{J_{<\gth}} f_\gg$.
 Finally ${\frak b} \subseteq {\frak b}^*$,
 because $\gr \in {\frak b} \implies f_\gg(\gr) \le f(\gr) < g(\gr)$.

 But now $\max(\pcf_W({\frak b}^*)) < \gth$,
 because otherwise we could find an ultrafilter $D$ on ${\frak b}^*$ with
 $\cf(\Pi {\frak b}^*/D) \ge \gth$. This contradicts ${\frak b}^* \in J_{<\gth}$.

  Applying the induction hypothesis, $ch^N_{{\frak b}^*} \in W$.
 Therefore 
\[
     ch^N_{\frak a} = f \restriction ({\frak a} \setminus {\frak b}^*) \cup ch^N_{{\frak b}^*} \in W,
\]  
 and we are done.
\end{proof}

   We conclude this section by proving that structures $N$ with
 most of the properties demanded in Lemma \ref{soft} are quite
 easily manufactured.

\begin{lemma} \label{manufacture}
 Let $W \subseteq V$ be two transitive class models of
 set theory. In $V$ let $\gk$ and $\gs$ be regular, with $\gk < \gs$.
 Assume that $\gk$-covering holds between $V$ and $W$, and
 that $\cf_W(\ga) \ge \gk^+_V \implies \cf_V(\ga) \ge \gk^+_V$ for all $\ga$. 
 Let  ${\cal M} = (H_\gs, \in, <^*, W \cap H_\gs)$ where $<^*$ is a
 wellordering of $H_\gs$.

 Let $\seq{N_i: i < \gb}$ be an increasing and continuous chain of
 substructures of ${\cal M}$ such that
\begin{enumerate}
\item $\gk \subseteq N_0$.
\item $\gb < \gk^+$ and  $\cf(\gb) = \gk$.
\item $\card{N_i} = \gk$ for all $i < \gb$.
\item $\seq{N_i : i \le j} \in N_{j+1}$ for all $j < \gb$.
\end{enumerate}

  Let $N =_{\rm def} \bigcup_{i < \gb} N_i$. Then 

\begin{enumerate}
\item If $\gth \in N \cap  REG^W$ with $\gth \ge \gk^+_V$, then
\begin{enumerate}
\item  $\sup(N \cap \gth) < \gth$. \label{z1}
\item $\cf_V(\sup(N \cap \gth)) = \gk$. \label{z2}
\end{enumerate}
\item For all ${\frak a} \subseteq REG^W \setminus \gk^+_V$ with ${\frak a} \in N \cap W$
 and $\card{{\frak a}} < \gk$, if $g \in \Pi {\frak a}$ and $g(\gr) < \sup(N \cap \gr)$
 for all $\gr \in {\frak a}$, then there is $h \in \Pi {\frak a} \cap N \cap W$ such that
 $g(\gr) < h(\gr)$ for all $\gr \in {\frak a}$. \label{z3}
\end{enumerate}
\end{lemma}

\begin{proof} We take each claim in turn.
 Since $\gth \ge \gk+_V$ and $\gth \in REG^W$, it follows from our assumptions
 that $\cf_V(\gth) > \gk$. $\card{N} = \gk$, so $\sup(N \cap \gth) < \gth$
 and we have proved claim \ref{z1}.
 Since $\gth \in N$, $\gth \in N_i$ for some $i < \gb$. For all $j > i$
 we have  $\gth, N_i \in N_j$, so that $\sup(N_i \cap \gth) \in N_j$
 and therefore $\sup(N_i \cap \gth) < \sup(N_j \cap \gth)$. 
 The sequence $\seq{\sup(N_j \cap \gth): i < j < \gb}$ is increasing
 and cofinal in $\sup(N \cap \gth)$ so
 $\cf(\sup(N \cap \gth)) = \cf(\gb) = \gk$. This proves claim \ref{z2}.

 Let ${\frak a}, g$ be as in claim \ref{z3} of the lemma. If ${\frak a} \in N_i$, a similar
 argument to that
 given in the last paragraph shows that for every $\gr \in {\frak a}$ the
 sequence $\seq{ \sup(N_j \cap \gr) : i < j < \gb}$ is increasing
 and cofinal in $\sup(N \cap \gr)$. Hence for every $\gr \in {\frak a}$ there
 is $j < \gb$ with $g(\gr) < \sup(N_j \cap \gr)$; since
  $\cf(\gb) = \gk > \card{{\frak a}}$ we may find a fixed $j$ such that
 $g(\gr) < \sup(N_j \cap \gr)$. Since $N_j \in N$, the function 
 $g^*: \gr \longmapsto \sup(N_j \cap \gr) \in N$. 

 Now we apply covering in a routine way
 to find a set $X \subseteq {\frak a} \times \bigcup {\frak a}$
 with $\card{X}_V < \gk$, $X \in W$ and $X \supseteq g^*$.
 Since $g^* \in N$ and
 $N \prec {\cal M}$ (a structure into which we built information
 about $W$) we may assume that $X \in N \cap W$. We define
 $h(\gr) = \sup \setof{\gb}{(\gr, \gb) \in X}$, and then
 $h \in N \cap W$ and also, since $\cf_V(\gr) > \gk$,
 $h \in \Pi {\frak a}$. Clearly $g(\gr) < g^*(\gr) \le h(\gr)$
 for all $\gr$, and we are done.
\end{proof}

  So the missing ingredient for applying Lemma \ref{soft} is the
 existence of a set $C$ with $C \in W$ and $C$ unbounded in $N \cap \gth$.
 In the next section we see two ways of guaranteeing the existence
 of such a $C$.

 Notice that the covering assumption in Lemma \ref{manufacture} can
 be weakened. All we need is that a set of size less than $\gk$ can
 be covered by a set of size less than $\gk^+$.

\section{Applying the pcf induction lemma}

We begin by showing that we can use Lemma \ref{soft} to prove
 instances of the principle $\rm{(J)'}$ which we now recall.

\medskip

\noindent ${\rm (J)'}_{\gk, \gl}$ :
For any ${\frak a} \subseteq \Reg^W \cap ( \gl \setminus
\gk^+)$ of cardinality $\gk$, if ${\frak a}$ is $W$-filtered, then for
every $f\in \Pi {\frak a}$ there exists $g \in \Pi {\frak a}$ such that $f \leq g$ and
$g$ is $W$-filtered.

\medskip

 As usual, $W$ and $V$ are transitive class
 models of ZFC with $W \subseteq V$, and $\gk$ is a regular cardinal
 in $V$. 

\begin{lemma} \label{applysoft} 
 Let ${\frak a} \subseteq REG^W \setminus (\gk+1)$ be a $W$-filtered set
 of size $\gk$, as witnessed by a
 filtration $\vec a$ such that $a_i \in W$ for all $i < \gk$. 
 Let $f \in \Pi {\frak a}$. Suppose that $f, {\frak a}, \vec a \in N$ where $N$ is a
 structure obeying the conclusion of Lemma \ref{soft},
 and such that $\max(\pcf_W(a_i)) \in N$ for all $i$. Then
 $ch^N_{\frak a}$ is $W$-filtered and $f \le ch^N_{\frak a}$.
\end{lemma}

\begin{proof} Since ${\frak a} \subseteq N$ (because $\card{{\frak a}} = \gk \subseteq N$)
 it is easy to see that $f \le ch^N_{\frak a}$. Since we built the filtration
 $\vec a$ into $N$ we see that for all $i < \gk$ we have $a_i \in N \cap W$,
 and we may apply Lemma \ref{soft} to conclude that $ch^N_{a_i} \in W$.
 Therefore $\seq{ch^N_{a_i}: i < \gk}$ gives a filtration of $ch^N_{\frak a}$
 and we are done.
\end{proof}

\begin{remark} The same idea could be used to derive the stronger property
 ${\rm (J)}$ defined in the introduction. That is, the assumption that
 $\vec a$ is increasing is never used.
\end{remark}

 Now we describe, in Lemmas \ref{app1} and \ref{app2},
 two ways of building structures $N$ that satisfy
 the hypotheses of Lemma \ref{soft}.

\begin{definition} Let $\gm$ and $\gth$ be regular cardinals with $\gm < \gth$.
 Then a {\em square sequence on $\gth$ for cofinalities less than $\gm$}
 is a sequence $\seq{C_\ga:\ga < \gth, \cf(\ga) < \gm}$  such that $C_\ga$ is
 club in $\ga$,
 $\ot{C_\ga} < \gm$ and $\gb \in \lim(C_\ga) \implies C_\gb = C_\ga \cap \gb$.
\end{definition}

\begin{lemma} \label{app1}
 Let $W$ and $V$ be two transitive class models of ZFC
 with $W \subseteq V$. Let $\gk$ be regular in $V$. Suppose that

\[
 W \models \hbox{``there is a square on $\gth$ for cofinalities 
 less than $\gk^+_V$''}
\]
 for every $W$-regular cardinal $\gth$. Suppose also
 $\cf_W(\ga) \ge \gk^+_V \implies \cf_V(\ga) \ge \gk^+_V$  
 for all $\ga$. 

  Let $\seq{N_i: i < \gb}$ be a sequence of substructures of
 $(H_\gs, \in, <^*, W \cap H_\gs)$ such that
\begin{enumerate}
\item $\cf(\gb) > \go$ and $\gb < \gk^+$.
\item $\gk \subseteq N_0$ and $\card{N_i} = \gk$.
\item $\seq{N_i: i \le j} \in N_{j+1}$ for all $j < \gb$.
\end{enumerate}
  Let $N = \bigcup_{i < \gb} N_i$, and let $\gth \in REG^W \cap N$ with
 $\gth \ge \gk^+_V$.
 If $\vec C \in W$ is a square on $\gth$ for cofinalities less than $\gk$
 then  $C_{\sup(N \cap \gth)} \subseteq N$.
\end{lemma}

\begin{proof} Let $\bar\gth =_{\rm def}  \sup(N \cap \gth)$, where
 $\bar\gth < \gth$ because $\card{N} = \gk$ and $\cf_V(\gth) > \gk$.
 $C_{\bar\gth}$ is defined because $\cf_W(\bar\gth) < \gk^+_V$.
 
 If $\gth \in N_i$ and we define $\gth_j = \sup(N_j \cap \gth)$
 then we may argue as in Lemma \ref{soft} that  $\seq{\gth_j: i \le j < \gb}$
 is increasing, continuous and cofinal in $\bar\gth$.
 Now since $\cf(\gb) > \go$ and $C_{\bar\gth}$ is club in
 $\bar\gth$ there is a club $D$ of $j < \gb$ such that
 $\gth_j \in \lim(C_{\bar\gth})$. For each such $j$,
 $C_{\gth_j} = C_{\bar\gth} \cap \gth_j$; since $\gth_j \in N$,
 $\card{C_{\gth_j}} \le \gk$ and $\gk \subseteq N$ we see that
 $C_{\gth_j} \subseteq N$, so
 $C_{\bar\gth} = \bigcup_{j \in D} C_{\gth_j}  \subseteq N$ and we are
 done.
\end{proof}

\begin{remark} By Lemma \ref{manufacture}, the structure $N$ defined in
 Lemma \ref{app1} obeys all the hypotheses of Lemma \ref{soft}.
\end{remark}

  We now describe another way of getting the hypotheses of Lemma
 \ref{soft} to hold. Here we drop the assumption of squares but
 pay for this by needing to assume more resemblance between $V$ and
 $W$.

\begin{lemma} \label{app2}
  Let $W$, $V$ be transitive class models of ZFC with $W \subseteq V$. Let $\gk$,
 $\gs$ be regular in $V$ with $\gk < \gs$ and let
 ${\cal M} = (H_\gs, \in, <^*, W \cap H_\gs)$ for some wellordering $<^*$ of $H_\gs$.
  Assume that $\gk^+_V = \gk^+_W$, $\gk^{++}_V = \gk^{++}_W$, 
 $\cf_W(\gt) \ge \gk^+ \implies \cf_V(\gt) \ge \gk^+$,
  $\cf_W(\gt) \ge \gk^{++} \implies \cf_V(\gt) \ge \gk^{++}$.

 Let $\seq{N_i: i < \gk^+}$ be a continuous and increasing chain of elementary
 submodels of ${\cal M}$ such that $\gk \subseteq N_0$,
 $\card{N_j} = \gk$, $\seq{N_i: i \le j} \in N_{j+1}$ for all $j$. 

 Then there exists $j$ such that for every $\gth \in N_j$ with $\gth \in REG^W$,
 $\gth \ge \gk^+$ there is $C \in W$ such that $C$ is unbounded in $N_j \cap \gth$.

\end{lemma}

\begin{proof}
   First we dismiss the case $\gth = \gk^+$. In this case there is not too
 much to prove because $N_j \cap \gth \in \gth$ for any $j$. So henceforth we
 assume that $\gth \ge \gk^{++}$.

  Next we claim that for a fixed $\gth \in N =_{\rm def} \bigcup_j N_j$,
 where $\gth \in REG^W \setminus \gk^{++}$,  there is a
 club $D$ in $\gk^+$ such that every $j \in D$ with $\cf(j) = \gk$ has the property
 claimed.

 Let $\gth \in N_i$, and fix $\seq{E_\gb: \gb < \gth} \in W \cap N_i$
 such that $E_\gb$ is club in $\gb$ and $\ot{E_\gb} = \cf_W(\gb)$. 
 Let $\gth_j =_{\rm def} \sup(\gth \cap N_j)$, so that as usual
 $\seq{\gth_j: i \le j < \gk^+}$ is continuous increasing and cofinal
 in $N \cap \gth$. 
 Our assumptions imply that $\cf_V(\gth) \ge \gk^{++}$, so that in particular
 $\bar\gth =_{\rm def} \sup(N \cap \gth) < \gth$ and $E_{\bar\gth}$ is defined.
 Our assumptions also imply that $\cf_W(\bar\gth) = \cf_V(\bar\gth) = \gk^+$,
 so that $\ot{E_{\bar\gth}} = \gk^+$.

 Now define $D \subseteq \gk^+$ by
\[
    D = \setof{j>i}{E_{\bar\gth} \cap \gth_j \cap N \subseteq N_j,
    \gth_j \in \lim(E_{\bar\gth})}.
\]
  It is easy to see that $D$ is club in $\gk^+$, the key point being that
 $\ot{E_{\bar\gth}} = \gk^+$ so that for each $j$ there is $k < \gk^+$
 with $E_{\bar\gth} \cap \gth_j \cap N \subseteq N_k$.

  Let $j \in D$ with $\cf(j) = \gk$, and let $C = E_{\gth_j} \cap E_{\bar\gth}$.
 Then $C \in W$ because $C$ is the intersection of two sets in $W$.
 By our assumptions again $\ot{E_{\gth_j}} = \cf_W(\gth_j) = \cf_V(\gth_j) = \gk$,
 so easily $\ot{C} = \gk$. Since $\vec E \in N_i$, $\gk \subseteq N_0$ and
 $\gth_j \in N_{j+1}$ we see that $E_{\gth_j} \subseteq N_{j+1} \subseteq N$.
 Now the key point is that
\[
   C = E_{\bar\gth} \cap E_{\gth_j} \subseteq E_{\bar\gth} \cap \gth_j \cap N \subseteq N_j,
\]
 as required.
 
   To finish the argument we just take a diagonal intersection. For each appropriate $\gth$
 fix $D_\gth$ club in $\gk^+$, such that every $j$ of cofinality $\gk$ in $D_\gth$
 is as desired.
 Now for each $i$ let $D^*_i = \bigcap_{\gth \in N_i} D_\gth$, and define a diagonal
 intersection
 $D^* = \setof{j < \gk^+}{i < j \implies j \in D^*_i}$.
 Then if $j \in D^*$ with $\cf(j) = \gk$,
 and $\gth \in N_j$, we see immediately that $\gth \in N_i$ for $i < j$ and thus
 $j \in D_\gth$.
\end{proof}

   We sum up the results of this section in a corollary.

\begin{corollary} \label{cor}
Let $W \subseteq V$ be two inner models of ZFC. Suppose that $\gk$, $\gs$ are
 regular cardinals in $V$ with $\gk < \gs$.
Suppose that $\cf_W(\ga) \ge \gk^+_V \implies \cf_V(\ga) \ge \gk^+_V$ for
 all $\ga < \gs$, and that either one of the following two assumptions holds:
\begin{enumerate}

\item \label{cor1} $W \models$ ``there is a square on $\gth$ for cofinalities less than $\gk^+_V$''
 for every $W$-regular $\gth$ such that $\gk^+_V \le \gth < \gs$.

\item \label{cor2} $\gk^+_V = \gk^+_W$, $\gk^{++}_V = \gk^{++}_W$, and
 $\cf_W(\ga) \ge \gk^{++} \implies \cf_V(\ga) \ge \gk^{++}$ for all
 $\ga < \gs$.
\end{enumerate}

  Then ${\rm (J)}_{\gk, \gs}$ holds.
\end{corollary}

\section{Conclusion}

 We can finally state the main  theorem.

\begin{theorem} \label{da_big_one}

  Let $W$ be an inner model of ZFC. Suppose that $\gk < \gl$ where $\gk$ is
 regular and $\gl$ is a cardinal in $W$. Suppose that
\begin{enumerate}
\item ${\rm (H)}_\gk$ holds. That is, every ordinal in $(\gk, \gk^+)$
 is $W$-filtered. 
\item There exists ${\cal P} \in W$ such that
 $W \models$ ``${\cal P} \subseteq [\gl]^{<\gk}$, $\card{{\cal P}} = \gl$''
and $V \models$ ``$\cal P$ is
 cofinal in $[\gl]^{<\gk}$''. \label{zz1}
\item ${\rm (I)}_{\gk, \gl}$ holds. That is, if
 $W \models$ ``${\frak a} \subseteq REG \cap \gl \setminus \gk^+$, $\card{\frak a} < \gk$''
 then $W \models ``\max(\pcf(\frak a)) \le \gl$''. 
\item One of the following holds:
\begin{enumerate} 
\item \label{ass1}
 ${\rm (D)_{\gk^+_V, \gl}}$ holds. That is, in $W$ there is a square sequence
 for points of cofinality less than $\gk^+_V$ on every regular cardinal in the interval
 $(\gk^+_V, \gl]$. Also 
 $\cf_W(\gt) \ge \gk^+ \implies \cf_V(\gt) \ge \gk^+$ for all $\gt < \gl$.

\item \label{ass2}
 $\gk^+_V = \gk^+_W$, $\gk^{++}_V = \gk^{++}_W$, and 
 $\cf_W(\gt) \ge \gk^+ \implies \cf_V(\gt) \ge \gk^+$,
 $\cf_W(\gt) \ge \gk^{++} \implies \cf_V(\gt) \ge \gk^{++}$
 for all $\gt < \gl$.
\end{enumerate}
\end{enumerate}

  Then strong $\gk$-covering holds between $V$ and $W$, for structures
 with underlying set $\gl$.

\end{theorem}

\begin{proof}
  The structure of this proof can be seen by looking at the picture in the
 introduction.
    Notice that assumption \ref{zz1} implies that ${\rm (B)}_{\gk, \gl}$ holds.
   By  Lemmas \ref{applysoft} and \ref{app1} (if we are assuming squares as in \ref{ass1}) 
 or \ref{applysoft} and \ref{app2} (if we are assuming correctness as in \ref{ass2})
 we have that ${\rm (J)'}_{\gk, \gl}$ holds. 
   By Lemma \ref{T21}, ${\rm (F)}_{\gl, \gk}$ holds. By Theorem \ref{FimpliesG},
 ${\rm (G)}_{\CP, \gk}$ holds. By Lemma \ref{Gimpliescovseq}, strong $\gk$-covering
 for sequences of subsets of $\gl$ holds. As in  the proof of Theorem \ref{painless},
 this implies that strong $\gk$-covering holds for structures on $\gl$.
\end{proof}

\section{Appendix on pcf}

  In this appendix we will prove the elementary facts about pcf theory
  used in the paper. For more information see the book \cite{CardinalArithmetic} or
  the survey paper \cite{BurkeMagidor}. 

\begin{definition} Let ${\frak a}$ be a set of regular cardinals such that
 $\card{{\frak a}}^+ < \min({\frak a})$. 
\begin{enumerate}
\item If $I$ is an ideal on ${\frak a}$ then
\begin{enumerate}
\item If  ${\frak b}, {\frak c} \subseteq {\frak a}$,
 ${\frak b} \subseteq_I {\frak c}$ iff ${\frak b} \setminus {\frak c} \in I$.
\item If $f, g \in \Pi {\frak a}$ then $f <_I g$ iff 
 $\setof{\gth \in {\frak a}}{g(\gth) \le f(\gth)} \in I$.
\end{enumerate}
\item If $F$ is a filter on $I$ and $f, g \in \Pi {\frak a}$ then
 $f <_F g$ iff $\setof{\gth}{f(\gth) < g(\gth)} \in \Pi {\frak a}$.
\item A strict partial ordering $(\FP, <_\FP)$ has {\em true cofinality}
 $\gl$ iff $\gl$ is regular and  there exists a sequence $\seq{p_i: i < \gl}$
 such that 
\begin{enumerate}
\item $i < j \implies p_i <_\FP p_j$.
\item $\forall p \in \FP \; \exists i < \gl \; p <_\FP p_i$.
\end{enumerate}
\item If $f, g \in \Pi {\frak a}$ then $f < g$ iff $f(\gth) < g(\gth)$ for
 all $\gth$.
\end{enumerate}
\end{definition}

  In general there is no guarantee that a partial ordering will
 have a true cofinality. If $\FP$ has a true cofinality it is
 easily seen to be unique, and we will write ``$\tcf(\FP)=\gl$''
 for the assertion that $\FP$ has true cofinality $\gl$;
 in the case that $\FP = (\Pi {\frak a}, <_I)$ for some ideal $I$ we
 will write ``$\tcf(\Pi {\frak a}/I) = \gl$''.
 The following lemma is easy, because an ultraproduct of cardinals
 will be linearly ordered and any linear ordering has a true cofinality. 

\begin{lemma} If $D$ is an ultrafilter on ${\frak a}$ then $(\Pi {\frak a}, <_D)$
 has a true cofinality.
\end{lemma}
 
  In this case we will write ``$\cf(\Pi {\frak a}/D) = \gl$''.

\begin{definition} Let ${\frak a}$ be a set of regular cardinals with 
 $\card{{\frak a}}^+ < \min({\frak a})$. Then
\begin{enumerate} 
\item $\pcf{{\frak a}}$ (Potential CoFinalities of ${\frak a}$) is the set of regular cardinals 
 $\gl$ such that $\cf(\Pi {\frak a}/D) = \gl$ for some ultrafilter $D$.
\item If $\gl$ is a cardinal (not necessarily regular) then
 $J_{<\gl}[{\frak a}]$ is the set of those ${\frak b} \subseteq {\frak a}$ such that
 ${\frak b} \in D \implies \cf(\Pi {\frak a}/D) < \gl$ for every ultrafilter $D$.
\end{enumerate}
\end{definition}

   Usually ${\frak a}$ will be clear from the context and we just write $J_{<\gl}$
 for $J_{<\gl}[{\frak a}]$. The following is the first key fact in pcf theory.

\begin{lemma} \label{key1}
 Let ${\frak a}$, $\gl$ be as above.
\begin{enumerate}
\item $J_{<\gl}$ is an ideal (possibly an improper one).
\item The poset $(\Pi {\frak a}, <_{J_{<\gl}})$ is {\em $\gl$-directed,}
 that is to say that if $F \subseteq \Pi {\frak a}$ and $\card{F} < \gl$
 there is $g \in \Pi {\frak a}$ such that $\forall f \in F \; f <_{J_{<\gl}} g$.
\end{enumerate}
\end{lemma}

\begin{proof} For brevity, let $J = J_{<\gl}$.
 The proof that $J$ is an ideal is fairly routine. For
 example let ${\frak b},{\frak c} \in J$ and ${\frak b} \cup {\frak c} \in D$ for some ultrafilter $D$.
 Because $D$ is an ultrafilter, either ${\frak b} \in D$ or ${\frak c} \in D$, and in either
 case $\cf(\Pi {\frak a}/D) < \gl$. 

 The proof of $\gl$-directedness goes by induction on $\card{F}$ for $F \subseteq \Pi {\frak a}$.
 Let $\card{F} = \gm$.
 If $\gm \le \card{{\frak a}}^+$ then we may define a bound $g(\gth) = \bigcup_{f \in F} (f(\gth)+1)$,
 so without loss of generality $\gm > \card{{\frak a}}^+$. If $\gm$ is singular then we may write
 $F = \bigcup_{i < \cf \gm} F_i$ with $\card{F_i} < \gm$, and then apply the induction
 hypothesis to find a bound $g_i$ for each $F_i$ and then a bound $g$ for all the $g_i$.
 So we may assume that $\card{{\frak a}}^+ < \gm = \cf(\gm) < \gl$.

 If $F$ is enumerated as $\seq{f_\ga: \ga < \gm}$, then we may use the induction
 hypothesis to  define inductively
 $f^*_\ga$ which is a $<_J$-upper bound for the set $\setof{f_\ga}{ \ga < \gm} \bigcup \setof{f^*_\ga}{ \ga < \gm}$.
 Then $\seq{f^*_\ga: \ga < \gm}$ is $<J$-increasing, and a bound for $\setof{f^*_\ga}{ \ga < \gm}$
 will be a bound for $\setof{f_\ga}{\ga<\gm}$. Relabeling, we may as well assume that we are trying
 to find an upper bound for a sequence $\seq{f_\ga: \ga < \gm}$ which is $<_J$-increasing.

 We will define a sequence of functions $g_\gb$ such that $\gb < \gg \implies g_\gb < g_\gg$,
 in such a way that for some $\gb < \card{{\frak a}}^+$ the function $g_\gb$ will be a bound for
 $\vec f$. $g_0(\gth) = 0$ for all $\gth$, and for limit $\gl < \card{{\frak a}}^+$ we define 
 $g_\gl = \sup_{\gb < \gl} g_\gb(\gth)$; $g_\gl \in \Pi {\frak a}$ because $\cf(\gth) = \gth > \card{{\frak a}}^+$
 for all $\gth \in {\frak a}$.

 Suppose that we have defined  $g_\gb$. If it fails to be a bound for the sequence $\vec f$ then
 $\setof{\gth}{f_{\ga(\gb)}(\gth) > g_\gb(\gth)} \in J^+$
 for some $\ga(\gb) < \gm$. 
 By the definition of $J$ we may choose
 $D$ an ultrafilter such that $g_\gb <_D f_{\ga(\gb)}$ and $\cf(\Pi {\frak a}/D) \ge \gl$
 (so necessarily $D \cap J = \emptyset$!)
 and we will then define $g_{\gb + 1}$ to be some function such that $g_{\gb +1} > g_\gb$
 and $g_{\gb+1}$ is an upper bound for $\vec f$ in $\Pi {\frak a}/D$. The key point here is
 that since $D \cap J = \emptyset$ and $\vec f$ is $<_J$-increasing we have 
 $g_\gb <_D f_\ga <_D g_{\gb+1}$ for all $\ga \ge \ga(\gb)$, so that in particular
 there exists $\gth$  with $g_\gb(\gth) < f_\ga(\gth) < g_{\gb+1}(\gth)$.

 Suppose for a contradiction that the construction of $\vec g$ runs for $\card{{\frak a}}^+$
 many steps. Choose $\ga \ge \sup_{\gb < \card{{\frak a}}^+} \ga(\gb)$. Then for all
 $\gb < \card{{\frak a}}^+$ there exists $\gth \in {\frak a}$ with $g_\gb(\gth) < f_\ga(\gth) < g_{\gb+1}(\gth)$.
 Some $\gth$ must occur twice, but this is absurd because $\vec g$ is a sequence which
 increases on every coordinate.
 
  We have shown that for some $\gb < \card{{\frak a}}^+$, $g_\gb$ is a bound. This
 concludes the proof.
\end{proof}
  
\begin{corollary} Let ${\frak a}$ be a set of regular cardinals with
 $\card{\frak a}^+ < \min(\frak a)$. Then
\begin{enumerate}
\item
 For every cardinal $\gl$, 
 $\cf(\Pi {\frak a}/D) < \gl$ iff $D \cap J_{<\gl} \neq \emptyset$.
\item $\pcf({\frak a})$ has a maximum element. 
\end{enumerate}
\end{corollary}

\begin{proof} The right to left direction in the first claim is immediate from the
 definition of $J_{<\gl}$, so suppose that 
 $D \cap J_{<\gl} = \emptyset$. In this case it follows immediately
 from the $\gl$-directedness of $(\Pi {\frak a}, <_{J_{<\gl}})$ that
 $(\Pi_{\frak a}, <_D)$ is also $\gl$-directed, so that clearly $\cf(\Pi {\frak a}/D) \ge \gl$.

 For the second claim suppose that $\pcf({\frak a})$ has no largest element, and let
 $\gm$ be the supremum of $\pcf({\frak a})$. ${\frak a} \notin J_{<\gl}$
 for each $\gl < \gm$, because ${\frak a} \in J_{<\gl} \implies \sup \pcf({\frak a}) \le \gl$.
 So $\bigcup_{\gl < \gm} J_{<\gl}$ is a proper ideal, and we may choose 
 $D$ an ultrafilter disjoint from it. But now $\cf(\Pi {\frak a}/D) = \gl$ for some
 $\gl < \gm$, and so $D \cap J_{<\gl^+} \neq \emptyset$, contradicting
 the choice of $D$.
\end{proof}

   Now we turn our attention to problems about true cofinalities.
 We need a technical lemma about unbounded sequences modulo some ideal
 $I$; intuitively the lemma says that an unbounded sequence can be split 
 into cofinal and bounded parts.

\begin{lemma} \label{key2}
 Let ${\frak a}$ be a set of regular cardinals, with $\card{{\frak a}}^+ < \min({\frak a})$.
 Let $I$ be an ideal on ${\frak a}$, let $\gm = \cf(\gm) > \card{{\frak a}}^+$, and
 $\seq{f_\ga: \ga < \gm}$ be increasing and unbounded in $(\Pi {\frak a}, <_I)$.
 Then there is a sequence $\seq{{\frak b}_\gg: \gg < \gm}$ of sets in $I^+$
 such that 
\begin{enumerate}
\item $\gg < \gd \implies {\frak b}_\gg \subseteq_I {\frak b}_\gd$.
\item $\seq{f_\ga \restriction {\frak b}_\gg: \ga < \gm}$ is cofinal in
 $\Pi {\frak b}_\gg/I$ for each $\gg$.
\item $\vec  f$ is bounded modulo the ideal generated by $I$ and
 the ${\frak b}_\gg$'s.
\end{enumerate}
\end{lemma}

\begin{proof} As in the proof of Lemma \ref{key1} we will build a sequence
 $\vec g$ of functions which are increasing on each coordinate, where
 $g_0 = 0$ and we take the pointwise supremum at each limit stage $\gb < \card{{\frak a}}^+$.
 Suppose we have defined $g_\gb$. Then since $\vec f$ is unbounded, the set
 ${\frak b}^\gb_\gg = \setof{\gth}{g_\gb(\gth) < f_\gg(\gth)}$ is in $I^+$ 
 for all $\gg$ sufficiently large (say $\gg \ge \gg(\gb)$). 
 Consider $\seq{{\frak b}^\gb_\gg: \gg \ge \gg(\gb)}$ as a candidate for the
 desired sequence $\vec {\frak b}$; clearly it is positive and increasing modulo
 $I$, and what is more $g_\gb$ will be a bound for $\vec f$ modulo the ideal
 generated by $I$ and $\seq{{\frak b}^\gb_\gg: \gg \ge \gg(\gb)}$. 

 So the construction is finished unless there is an $\gg^*(\gb) \ge \gg(\gb)$ such that
 $\seq{f_\ga \restriction {\frak b}^\gb_{\gg^*(\gb)}: \ga < \gm}$ fails to be cofinal.
 In this case we will choose $g_{\gb+1} > g_\gb$ to be a witness
 to this failure of cofinalness, which is to say that
 $\setof{\gth \in {\frak b}^\gb_{\gg^*(\gb)}}{g_{\gb+1}(\gth) > f_\ga(\gth)} \in I^+$
 for all $\ga$.
 The key point is that (since the ${\frak b}^\gb_\gg$ are increasing modulo $I$
 with $\gg$) for all $\gg \ge \gg^*(\gb)$ and all $\ga$ we have 
 $\setof{\gth \in {\frak b}^\gb_\gg}{g_{\gb+1}(\gth) > f_\ga(\gth)} \in I^+$.

 Now suppose that the construction runs for $\card{{\frak a}}^+$ many steps. 
 Choose $\gg \ge \sup_{\gb < \card{{\frak a}}^+} \gg^*(\gb)$. Applying the conclusion
 of the last paragraph for $\ga =\gg$,  
  $\setof{\gth \in {\frak b}^\gb_\gg}{g_{\gb+1}(\gth) > f_\gg(\gth)} \in I^+$,
 so in particular there is $\gth \in {\frak a}$ such that
 $g_\gb(\gth) < f_\gg(\gth) < g_{\gb+1}(\gth)$. This leads to a contradiction
 exactly as in the proof of Lemma \ref{key1}.

 This shows that at some stage $\gb < \card{{\frak a}}^+$ the construction
 terminates, giving a sequence $\vec {\frak b}$ as desired.
\end{proof}

   Using the facts above we can derive the key fact about pcf which is being
 used in this paper.

\begin{theorem} \label{bop} Let $\gth = \max\pcf({\frak a})$. Then
 $\tcf(\Pi {\frak a}/J_{<\gth}) = \gth$. 
\end{theorem}

\begin{proof} Define  $J = J_{<\gth} \bigcup
                                      \setof{{\frak b} \notin J_{<\gth}}{\tcf(\Pi {\frak b}/J_{<\gth}) = \gth}$.
 It is easy to check that $J$ is a (possibly improper) ideal.
 If ${\frak a} \in J$ we are done, otherwise let us choose $D$ an ultrafilter on ${\frak a}$ with
 $J \cap D =\emptyset$. Since $J_{<\gth} \cap D = \emptyset$ and $\gth = \max\pcf({\frak a})$
 it is easy to see that $\cf(\Pi {\frak a}/D) = \gth$.

 Let us choose $\seq{f_\ga: \ga < \gth}$ which is increasing and cofinal in $(\Pi {\frak a}, <_D)$.
 Since $D \cap J_{<\gth} = \emptyset$ and $\Pi {\frak a}/J_{<\gth}$ is $\gth$-directed we may
 also assume that $\vec f$ is $<J_{<\gth}$-increasing. It is easy to see that since
 $\vec f$ is cofinal modulo $D$ it is unbounded modulo $J_{<\gth}$, so that we may
 apply Lemma \ref{key2}.

 Let $\vec {\frak b}$ be the sequence of sets given by Lemma \ref{key2}. If ${\frak b}_\gg \notin D$
 for all $\gg$ then the ideal generated by $J_{<\gth}$ and $\vec {\frak b}$ is contained
 in the ideal dual to $D$, contradicting the choice of $\vec f$ to be cofinal modulo
 $D$. If on the other hand there is a $\gg$ with ${\frak b}_\gg \in D$, then 
 $\seq{f_\ga \restriction {\frak b}_\gg: \ga < \gth}$ witnesses that
 $\tcf(\Pi {\frak b}_\gg/I) = J_\gth$, so ${\frak b}_\gg \in D \cap J$ contradicting the
 choice of $D$ disjoint from $J$.

 This contradiction shows that ${\frak a} \in J$, hence $\tcf(\Pi {\frak a}/J_{<\gth}) = \gth$ and
 we are done.
\end{proof}

 \end{document}